\documentclass[11pt]{amsart}
\usepackage{amsfonts,amssymb,amsmath,amsthm,here}

\usepackage{color}  
\usepackage{array}
\usepackage{graphicx}
\usepackage{float} 
\usepackage{caption}
\newfloat{graph}{H}{log}
\floatname{graph}{Graph}

\theoremstyle{plain}
\newtheorem{theo}{Theorem}

\newtheorem{lem}{Lemma}
\newtheorem{prop}{Proposition}

\theoremstyle{definition}

\newtheorem{conj}{Conjecture}

\makeatletter

\@addtoreset{equation}{section}
\makeatother

\makeatletter
\@namedef{subjclassname@2020}{
  \textup{2020} Mathematics Subject Classification}
\makeatother

\author{Miho Aoki}
\address{Department of Mathematics,
Interdisciplinary Faculty of Science and Engineering,
Shimane University,
Matsue, Shimane, 690-8504, Japan}
\email{aoki@riko.shimane-u.ac.jp}

\thanks{This work was supported by JSPS KAKENHI Grant Number JP26K06765.} 

\title[Distribution of Primes in Galois Extensions of Global Fields]{Distribution of Primes in Galois Extensions of Global Fields with Explicit Constant Terms}
\date{}

\begin{document}
\maketitle

\begin{abstract}
We study the distribution of primes among Frobenius conjugacy classes in finite Galois extensions of global fields. In earlier work, the author and Koyama proved that
the existence and nonvanishing of the limit of partial Euler products on the critical line for Artin \(L\)-functions
(the Deep Riemann Hypothesis (A), or DRH(A)) are equivalent to an asymptotic formula for weighted prime counting functions with an unspecified constant term. 
In this article, we prove that a refined asymptotic formula is equivalent to the full version of DRH, namely the convergence of the Euler products on the critical line (DRH(A)) and the determination of their limiting values (DRH(B)). Here, the constant term is given explicitly in terms of the residue of the zeta function, central derivatives of Artin \(L\)-functions, ramified primes, and higher prime-power contributions.
\end{abstract}
\section{Introduction}

In a finite Galois extension of global fields, the Chebotarev density theorem guarantees that the prime ideals whose Frobenius conjugacy classes fall into a given conjugacy class appear asymptotically in proportion to the size of that class. However, up to any finite norm bound, the actual counts of prime ideals corresponding to different conjugacy classes of the same size frequently exhibit persistent disparities. The most classical manifestation of this phenomenon is Chebyshev's bias~\cite{Ch}: while the prime number theorem for arithmetic progressions (a special case of the Chebotarev density theorem) ensures that primes congruent to \(3\) and \(1\) modulo \(4\) each account for asymptotically half of all primes as the upper bound grows without limit, primes congruent to 3 modulo 4 often outnumber those congruent to \(1\) modulo \(4\) up to a given finite threshold.
Rubinstein and Sarnak \cite{RS} established a rigorous framework to quantify this bias in terms of logarithmic density.
Assuming the Generalized Riemann Hypothesis and the linear independence of the zeros of Dirichlet \(L\)-functions, Rubinstein and Sarnak showed that primes congruent to \(3\) modulo \(4\) outnumber those congruent to \(1\) modulo \(4\) across roughly \(99.59\) percent of the logarithmic scale.

The author and Koyama \cite{AK} formulated Chebyshev's bias using weighted prime counting functions and obtained asymptotic formulas under the convergence of the partial Euler products on the critical line, namely DRH(A). We established a criterion for the existence of a bias among Frobenius elements in Galois extensions of global fields. As applications, we showed that in abelian extensions, there is a bias against splitting primes and principal prime ideals. 
In the positive characteristic case, DRH is known to hold unconditionally, having been essentially established by Conrad~\cite{C} (see \cite{KKK} for the proof in full generality for $\mathrm{GL}_n$). 
Numerous results under the assumption of DRH have been obtained by Koyama, Kurokawa, Okumura, Sheth, and others \cite{K, KK0, KK1, KS, KS1, O}.
In the case of elliptic curves over function fields, Kaneko and Koyama~\cite{KK} applied this weighted approach to show unconditionally that the direction of the bias is directly dictated by the algebraic rank of the curve. Pushing this paradigm to higher-dimensional geometric objects, Bucur, Kedlaya, and Sheth~\cite{BKS} recently formulated a Birch--Swinnerton-Dyer-type asymptotic for products of point counts of higher-genus curves over finite fields.

In this article, we provide a refined asymptotic formula that is equivalent to the full version of DRH, namely the convergence of the Euler product (DRH(A)) and its limiting value (DRH(B)). 
Let \(K\) be a global field and  let \(S_{K}\)  denote the set of all finite primes (non-zero prime ideals) if \(K\) is an algebraic number field, and the set of all primes if \(K\) is a global function field.
In the case of a global function field, we assume that \(K\) is a finite extension of \(\mathbb{F}_q(T)\), and that the finite field \(\mathbb{F}_{q}\) is algebraically closed in \(K\). That is, the equality \(\{a \in K \mid a \text{ is algebraic over } \mathbb{F}_q\} = \mathbb{F}_q\) holds.
 For $\mathfrak p\in S_K$, define
\[
N\mathfrak p=
\begin{cases}
|O_K/\mathfrak p|, & \text{if \(K\) is a number field},\\
q^{\deg(\mathfrak p)}, & \text{if \(K\) is a global function field},
\end{cases}
\]
where \(O_K\) is the ring of integers of \(K\).
Let
\[
\rho:\operatorname{Gal}(K^{\mathrm{sep}}/K)\longrightarrow \operatorname{Aut}_{\mathbb C}(V)
\]
be a nontrivial irreducible Artin representation. Its Artin $L$-function is

\[
L_K(s,\rho)
=
\prod_{\mathfrak p\in S_K}
\det\left(I-N\mathfrak p^{-s}\rho(\operatorname{Frob}_{\mathfrak p})\mid  V^{I_{\mathfrak p}}  \right)^{-1},
\qquad (\operatorname{Re}(s)>1).
\]
Put
\[
m_\rho=\operatorname{ord}_{s=1/2}L_K(s,\rho),
\]
and define
\[
\nu (\rho)
=
\operatorname{mult}(\mathbf{1},\operatorname{Sym}^2\rho)
-
\operatorname{mult}(\mathbf{1},\wedge^2\rho),
\]
where \(\operatorname{mult}(\mathbf{1},\rho')\) is the multiplicity of the trivial representation \( \mathbf{1} \) in $\rho'$.
The Deep Riemann Hypothesis (DRH) \cite{KKK} consists of two related statements: DRH(A) asserts the existence and nonvanishing of the Euler product on the critical line, while DRH(B) identifies its limiting value.
\begin{conj}[Deep Riemann Hypothesis (DRH)]\label{DRH}
Put $m=m_\rho$.
Then the limit 
\begin{equation}\label{limit}
\lim_{x \to \infty} \Bigg((\log x)^{m} \prod_{N\mathfrak p \le x} 
\det\left(1-\rho(\mathrm{Frob}_{\mathfrak p})N\mathfrak p^{-1/2} | V^{I_{\mathfrak p}}  \right)^{-1}  \Bigg)
\end{equation}
satisfies the following:
\begin{description}
\item[DRH(A)] The limit \eqref{limit} exists and is nonzero.
\item[DRH(B)] The limit \eqref{limit} satisfies the following identity:
\[
\lim_{x \to \infty} \Bigg((\log x)^{m} \prod_{N\mathfrak p \le x} 
\det\left(1-\rho(\mathrm{Frob}_{\mathfrak p})N\mathfrak p^{-1/2} | V^{I_{\mathfrak p}} \right)^{-1} \Bigg)
 = \frac{\sqrt{2}^{\nu(\rho)} L_{K}^{(m)}(1/2, \rho)}{e^{m \gamma} m!},
\]
where $L_K^{(m_\rho)}(1/2,\rho)$ denotes the first nonzero derivative at $s=1/2$.
\end{description}
\end{conj}
%

Conrad \cite{C} showed that \(\mathrm{DRH}(\mathrm{A})\) implies the Generalized Riemann Hypothesis (\(\mathrm{GRH}\)) 
and that the convergence of the partial Euler products at \(s=1/2\) is equivalent to their convergence at every \(s\) on the line \(\mathrm{Re}(s)=1/2\).
In the case of global function fields, both \(\mathrm{DRH}(\mathrm{A})\) and \(\mathrm{DRH}(\mathrm{B})\) are known to hold unconditionally, thanks to the work of Conrad \cite{C} via the Weil conjectures proved by Deligne \cite{D} and further explicit proofs by Kaneko, Koyama, and Kurokawa \cite{KKK}. 
 Furthermore, Kimura, Koyama, and Kurokawa \cite{KKK1} provided explicit numerical examples verifying \(\mathrm{DRH}\) for quadratic Dirichlet characters.
For the Riemann zeta-function, the relationship between the convergence of partial Euler products on the critical line and prime distribution was rigorously established  by Akatsuka \cite{A}, who treated the pole at $s=1$ via an explicit regularization.

Let \(L/K\) be a finite Galois extension, where we assume that the finite field \(\mathbb{F}_{q}\) is algebraically closed in \(L\) if \(K\) is a global function field and put
\(G=\operatorname{Gal}(L/K).\)
For \(\sigma\in G\), let \(c_\sigma\) denote its conjugacy class, and let \(S_\sigma\) be the set of primes \(\mathfrak p\) of \(K\) unramified in \(L\)  whose Frobenius conjugacy class equals \(c_\sigma\). Chebotarev's density theorem implies that \(S_\sigma\) has density
\(
|c_\sigma|/[L:K]
\).
Following \cite{AK}, we define the weighted prime counting functions as
\[
\pi_{\frac{1}{2},K}(x)
=
\sum_{\substack{\mathfrak p\in S_K \atop N\mathfrak p\leq x}}
\frac{1}{\sqrt{N\mathfrak p}},
\]
and for \(\sigma \in G\)
\[
\pi_{\frac{1}{2}}(x;\sigma)
=
\sum_{\substack{\mathfrak p\in S_\sigma \atop N\mathfrak p\leq x}}
\frac{1}{\sqrt{N\mathfrak p}}.
\]
For each \(\sigma\in G\), define
\[
M(\sigma)
=
\frac{1}{2}
\sum_{\substack{\rho\ne \mathbf{1} \\ {\rm irred.}}}
\chi_\rho(\sigma)
\nu(\rho), \quad m(\sigma)
=
\sum_{\substack{\rho\ne \mathbf{1} \\ {\rm irred.}} }
\chi_\rho(\sigma)m_\rho,
\]
where $\chi_{\rho}$ is  the character of $\rho$ and
$\rho$ runs over all nontrivial irreducible representations of \(\mathrm{Gal}(L/K)\). 
The quantity $M(\sigma)$ is defined algebraically, whereas $m(\sigma)$ is defined analytically since it is determined by orders of vanishing of Artin $L$-functions at the central point.
We can easily check that $M(\sigma), m(\sigma ) \in \mathbb R$.

The weighted prime counting functions introduced above play a pivotal role in the study of Chebyshev's bias. 
The author and Koyama \cite{AK} derived an explicit asymptotic formula for the weighted prime counting function, proving that its validity is equivalent to the assertion of DRH(A).

\begin{theo}[Aoki--Koyama]\label{theo:AK}
Let \(L/K\) be a finite Galois extension of global fields  with Galois group \(G\). The following are equivalent:
\begin{enumerate}
\item DRH(A) holds for every nontrivial irreducible representation of \(G\).
\item For every \(\sigma\in G\), there exists a constant \(C(\sigma)\) such that
\[
\pi_{\frac{1}{2},K}(x) -\frac{[L:K]}{|c_{\sigma}|} \pi_{\frac{1}{2}} (x;\sigma)
=
\bigl(M(\sigma)+m(\sigma)\bigr)\log\log x
+C(\sigma)+o(1)
\qquad (x\to\infty).
\]
\end{enumerate}
\end{theo}
This theorem identifies the coefficient of \(\log\log x\) as the quantitative measure of weighted Chebyshev bias. In particular,
\[
M(\sigma)+m(\sigma)=0
\]
means that the Frobenius class \(c_{\sigma} \) has no bias. If this coefficient is positive, primes with Frobenius class \(c_{\sigma}\) occur relatively late; if it is negative, they occur relatively early.

\section{Explicit Constants and DRH(B)}

We  refine the framework of Theorem~\ref{theo:AK} by incorporating the exact limiting value of the Euler products (DRH(B)), which allows us to fully formulate Chebyshev's bias with the explicit constant term \(C(\sigma)\). We show that the validity of this refined asymptotic formula is equivalent to both DRH(A) and DRH(B).
Let 
\[
\zeta_K(s) =\prod_{\mathfrak{p} \in S_K} \left( 1-\frac{1}{N\mathfrak p^s} \right)^{-1} \quad ({\mathrm Re}(s) >1)
\]
be the zeta function of \(K\), and let

\[
\alpha_K=\operatorname*{Res}_{s=1}\zeta_K(s) =\begin{cases} 
\dfrac{2^{r_1} (2\pi )^{r_2 }h_K R_K}{w_K  \sqrt{|D_K|}}  & \text{(if \(K \) is a number field)},  \\
& \\
\dfrac{h_K}{q^g (1-q^{-1})\log{q}}  & \text{(if \(K\) is a global function field)},
\end{cases}
\]
where \(r_1\) and \(r_2\) denote the number of real embeddings and pairs of complex embeddings of \(K\), respectively; 
\(h_K\) is the class number if \(K\) is a number field, and the order of the group of divisor classes of degree zero if \(K\) is a global function field;
\(R_K\) is the regulator; \(w_K\) is the number of roots of unity in \(K\);  \(D_K\) is the absolute discriminant; and \(g\) is the genus of the global function field \(K\).

Define
\[
c_K
=-\sum_{\mathfrak p\in S_K}
\left(
\frac{1}{N\mathfrak p}
+\log\left(1-\frac{1}{N\mathfrak p}\right)
\right).
\]
The infinite sum on the right-hand side converges (see \cite[Lemma 2.5]{R} for the case of algebraic number fields, and \cite[Lemma 1.2]{R} for global function fields).
Let $S_{\mathrm{ram}}$ be the set of primes of $K$  ramified in $L/K$, and put
\[
R_{L/K}=\sum_{\mathfrak p\in S_{\mathrm{ram}}}
\frac{1}{\sqrt{N\mathfrak p}}.
\]
\begin{theo}\label{theo:main}
Let \(L/K\) be a finite Galois extension of global fields with Galois group \(G\). The following are equivalent:
\begin{enumerate}
\item DRH(A) and DRH(B) hold for every nontrivial irreducible representation of \(G\).
\item For every \(\sigma\in G\), we have
\[
\pi_{\frac{1}{2},K}(x) -\frac{[L:K]}{|c_{\sigma}|} \pi_{\frac{1}{2}} (x;\sigma)
=
\bigl(M(\sigma)+m(\sigma)\bigr)\log\log x
+C(\sigma)+o(1) \qquad (x\to \infty),
\]
with
\begin{align*}
C(\sigma)
&= \bigl(M(\sigma)+m(\sigma)\bigr)\gamma+R_{L/K}-M(\sigma)\bigl(\log 2+c_K-\log\alpha_K\bigr) \\
&\quad -\sum_{\substack{\rho \neq \mathbf{1} \\ \mathrm{irred.}} }
\overline{ \chi_\rho (\sigma) }
\left(
\log\frac{L_K^{(m)}(1/2,\rho)}{m!}
-c(\rho)
\right),
\end{align*}
where \(\gamma\) is Euler's constant, \(m=m_\rho\), and \(c(\rho)\) is a constant arising from ramified primes and higher prime-power contributions. This constant is explicitly given by 
\begin{align*}
c(\rho) &= \sum_{\mathfrak{p} \in S_{\mathrm{ram}}} \sum_{\tau \in I_{\mathfrak{p}}(L/K)} \frac{\chi_{\rho} (\tau \operatorname{Frob}_{\mathfrak{p}})}{|I_{\mathfrak{p}}(L/K)| N\mathfrak{p}^{1/2}} \\
&\quad +
\frac{1}{2} \left\{ \sum_{\substack{\rho' \neq \mathbf{1} \\ \mathrm{irred.}}} \left( \bigl(\operatorname{mult}(\rho', \operatorname{Sym}^2\rho)- \operatorname{mult}(\rho', \wedge^2\rho)\bigr) \sum_{\mathfrak{p} \in S_K\setminus S_{\mathrm{ram}}} \frac{\chi_{\rho'} (\operatorname{Frob}_{\mathfrak{p}})}{N\mathfrak{p}} \right) \right. \\
&\quad \left. + \sum_{\mathfrak{p} \in S_{\mathrm{ram}}}
\frac{1}{N\mathfrak{p}} \left( \sum_{\tau \in I_{\mathfrak{p}}(L/K)}\frac{\chi_{\rho} (\tau \operatorname{Frob}_{\mathfrak{p}}^2)}{|I_{\mathfrak{p}}(L/K)|} -\nu(\rho) \right)\right\} \\
&\quad + \sum_{\mathfrak{p} \in S_K} \sum_{\tau \in I_{\mathfrak{p}}(L/K)} \sum_{k=3}^{\infty} \frac{\chi_{\rho} (\tau \operatorname{Frob}_{\mathfrak{p}}^k)}{k |I_{\mathfrak{p}}(L/K)| N\mathfrak{p}^{k/2}},
\end{align*}
where \(I_{\mathfrak{p}}(L/K)\) denotes the inertia group of a prime \(\mathfrak{P}\) of \(L\) above \(\mathfrak{p}\).
\end{enumerate}
\end{theo}

In the lemma stated below, we verify that the infinite series contained in the expression for \(c(\rho)\) in Theorem~\ref{theo:main} converge.
\begin{lem}\label{lem:const}
Under the notation and settings of Theorem~\ref{theo:main},
for any nontrivial irreducible representation \(\rho\) of \(G\),
the two infinite series contained in the expression for \(c(\rho)\) in Theorem~\ref{theo:main}:
\[
\sum_{\mathfrak{p} \in S_K\setminus S_{\mathrm{ram}}} \frac{\chi_{\rho} (\operatorname{Frob}_{\mathfrak{p}})}{N\mathfrak{p}} 
\]
and
\[
\sum_{\mathfrak{p} \in S_K} \sum_{\tau \in I_{\mathfrak{p}}(L/K)} \sum_{k=3}^{\infty}\frac{\chi_{\rho} (\tau \operatorname{Frob}_{\mathfrak{p}}^k)}{k |I_{\mathfrak{p}}(L/K)| N\mathfrak{p}^{k/2}}
\]
converge.
\end{lem}
\begin{proof}
First, we show that the first series converges.
For any nontrivial irreducible representation \(\rho\), the Artin \(L\)-function \(L_K(s,\rho)\) converges to a nonzero value at \(s = 1\) (for number fields see \cite[Cor.~5.47]{IK}; for global function fields this follows from the Weil conjectures  proved by Deligne \cite{D}). Hence, 
\[
c = \log\left( \prod_{\mathfrak{p} \in S_K} \det\bigl(I-N\mathfrak{p}^{-1} \rho(\operatorname{Frob}_{\mathfrak{p}}) \mid V^{I_{\mathfrak{p}}}  \bigr)^{-1} \right)
\]
is a constant. Taking the Taylor expansion at each unramified place, we obtain the following formula:
\begin{equation}\label{eq:c}
c = \sum_{\mathfrak{p} \in S_K \setminus S_{\mathrm{ram}}} \sum_{k=1}^{\infty}\frac{\chi_{\rho} (\operatorname{Frob}_{\mathfrak{p}}^k)}{kN\mathfrak{p}^k} + 
\log\left( \prod_{\mathfrak{p} \in S_{\mathrm{ram}}} \det\bigl(I-N\mathfrak{p}^{-1} \rho(\operatorname{Frob}_{\mathfrak{p}}) \mid V^{I_{\mathfrak{p}}} \bigr)^{-1} \right).
\end{equation}
Since
\begin{align*}
\left| \sum_{\mathfrak{p} \in S_K \setminus S_{\mathrm{ram}}} \sum_{k=2}^{\infty} \frac{\chi_{\rho} (\operatorname{Frob}_{\mathfrak{p}}^k)}{kN\mathfrak{p}^k} \right| 
&\leq \sum_{\mathfrak{p} \in S_K \setminus S_{\mathrm{ram}}} \sum_{k=2}^{\infty} \left| \frac{\chi_{\rho} (\operatorname{Frob}_{\mathfrak{p}}^k)}{kN\mathfrak{p}^k} \right| 
\leq \sum_{\mathfrak{p} \in S_K \setminus S_{\mathrm{ram}}} \sum_{k=2}^{\infty} \frac{\deg \rho}{kN\mathfrak{p}^k} \\
&\leq \frac{\deg \rho}{2} \sum_{\mathfrak{p} \in S_K \setminus S_{\mathrm{ram}}} \sum_{k=2}^{\infty} \frac{1}{N\mathfrak{p}^k} 
\leq \deg \rho \times \sum_{\mathfrak{p} \in S_K \setminus S_{\mathrm{ram}}} \frac{1}{N\mathfrak{p}^2} < \deg \rho \times \zeta_K(2) < \infty,
\end{align*}
it follows from \eqref{eq:c} that the series 
\[
\sum_{\mathfrak{p} \in S_K\setminus S_{\mathrm{ram}}} \frac{\chi_{\rho} (\operatorname{Frob}_{\mathfrak{p}})}{N\mathfrak{p}} 
\]
converges.

Next, the convergence of the second series follows from the following inequalities:
\begin{align*}
\left| \sum_{\mathfrak{p} \in S_K} \sum_{\tau \in I_{\mathfrak{p}}} \sum_{k=3}^{\infty} \frac{\chi_{\rho} (\tau \operatorname{Frob}_{\mathfrak{p}}^k)}{k |I_{\mathfrak{p}}| N\mathfrak{p}^{k/2}} \right| 
&\leq \sum_{\mathfrak{p} \in S_K} \sum_{\tau \in I_{\mathfrak{p}}} \sum_{k=3}^{\infty} \left| \frac{\chi_{\rho} (\tau \operatorname{Frob}_{\mathfrak{p}}^k)}{k |I_{\mathfrak{p}}| N\mathfrak{p}^{k/2}} \right| 
\leq \sum_{\mathfrak{p} \in S_K} \sum_{k=3}^{\infty} \frac{\deg \rho}{k N\mathfrak{p}^{k/2}} \\
&\leq \frac{\deg \rho}{3} \sum_{\mathfrak{p} \in S_K} \sum_{k=3}^{\infty} \frac{1}{N\mathfrak{p}^{k/2}} 
\leq \frac{4 \deg \rho}{3} \sum_{\mathfrak{p} \in S_K} \frac{1}{N\mathfrak{p}^{3/2}} < \frac{4 \deg \rho}{3} \zeta_K(3/2) < \infty,
\end{align*}
where \(I_{\mathfrak{p}} = I_{\mathfrak{p}}(L/K)\).
\end{proof}
\begin{lem}\label{lem:const-sum}
Under the notation and settings of Theorem~\ref{theo:main},
the constants \(C(\sigma)\) satisfy the identity:
\[
\sum_{\sigma\in G} C(\sigma) = [L:K]R_{L/K}.
\]
\end{lem}
\begin{proof}
Since \(\sum_{\sigma \in G} \chi_{\rho}(\sigma) = 0\) for any nontrivial irreducible representation \(\rho\), we obtain \(\sum_{\sigma \in G} M(\sigma) = \sum_{\sigma \in G} m(\sigma) = 0\). Combining these equalities with the expression for \(C(\sigma)\) in Theorem~\ref{theo:main} yields the claim.
\end{proof}

\section{Proof of Theorem}

In this section, we give the proof of Theorem~\ref{theo:main}.
We use the following result, which Rosen used to establish a generalized Mertens' theorem (see \cite[Lemmas~2.4 and 2.5]{R} for the number field case, 
and \cite[Lemmas~1.1 and 1.2]{R} for the global function field case).
\begin{theo}[Rosen]\label{theo:Rosen}
Let \(K\) be a global field. We have 
\[
\sum_{N\mathfrak{p} \leq x}\frac{1}{N\mathfrak{p}} = \log\log x + M_K + O\left(\frac{1}{\log x} \right) \qquad (x\to \infty),
\]
where 
\[
M_K = \sum_{\mathfrak{p} \in S_K} \left(\frac{1}{N\mathfrak{p}} + \log \left(1-\frac{1}{N\mathfrak{p}} \right) \right) + \gamma + \log\alpha_K 
\]
is the Mertens constant for \(K\).
\end{theo}
\begin{prop}\label{prop:main}
Let \(L/K\) be a finite Galois extension of global fields with Galois group \(G\), and let \(\rho\) be a nontrivial irreducible representation of \(G\).
Put \(m=m_{\rho}\).
Then DRH(A) and DRH(B) hold for \(\rho\) if and only if the following asymptotic formula holds:
\begin{align*}
\sum_{\substack{N\mathfrak{p}\leq x \\ \mathfrak{p} \in S_K \setminus S_{\mathrm{ram}}}} \frac{\chi_{\rho}(\operatorname{Frob}_{\mathfrak{p}})}{N\mathfrak{p}^{1/2}} 
&= -\left(\frac{\nu(\rho)}{2} + m \right) (\log\log x + \gamma) + \frac{\nu(\rho)}{2} (\log 2 + c_K - \log\alpha_K) \\
&\quad + \log \frac{L_K^{(m)}(1/2,\rho)}{m!} - c(\rho) + o(1) \qquad (x\to \infty).
\end{align*}
\end{prop}

By a well-known property of the trace, we have the following lemma.
\begin{lem}\label{lem:Tr}
Under the notation and settings of Proposition~\ref{prop:main}, for any \(\mathfrak{p}\in S_K\) and any \(k\in \mathbb{Z}_{\geq 1}\), we have
\[
\operatorname{Tr}\bigl(\rho(\operatorname{Frob}_{\mathfrak{p}}^k) \big|_{V^{I_{\mathfrak{p}}}}\bigr) = \frac{1}{|I_{\mathfrak{p}}|} \sum_{\tau \in I_{\mathfrak{p}}} \chi_{\rho}(\tau \operatorname{Frob}_{\mathfrak{p}}^k),
\]
where \(I_{\mathfrak{p}} = I_{\mathfrak{p}}(L/K)\).
\end{lem}
\begin{proof}
Since the linear map \(\rho(\operatorname{Frob}_{\mathfrak{p}}^k) : V \to V\) 
satisfies \(\rho(\operatorname{Frob}_{\mathfrak{p}}^k) v \in V^{I_{\mathfrak{p}}}\) 
for any \(v \in V^{I_{\mathfrak{p}}}\), it induces a linear map 
\(\rho(\operatorname{Frob}_{\mathfrak{p}}^k)\big|_{V^{I_{\mathfrak{p}}}} : V^{I_{\mathfrak{p}}} \to V^{I_{\mathfrak{p}}}\). Put 
\[
P = \frac{1}{|I_{\mathfrak{p}}|} \sum_{\tau \in I_{\mathfrak{p}}} \rho(\tau).
\]
Then we have \(P\big|_{V^{I_{\mathfrak{p}}}} = \mathrm{id}_{V^{I_{\mathfrak{p}}}}\) and \(P(V) = V^{I_{\mathfrak{p}}}\).
Hence, by standard properties of the trace, we obtain
\begin{align*}
\operatorname{Tr}\bigl(\rho(\operatorname{Frob}_{\mathfrak{p}}^k) \big|_{V^{I_{\mathfrak{p}}}}\bigr) 
&= \operatorname{Tr}\bigl(\rho(\operatorname{Frob}_{\mathfrak{p}}^k)\circ P\bigr) 
= \operatorname{Tr}\bigl(P \circ \rho(\operatorname{Frob}_{\mathfrak{p}}^k)\bigr) \\
&= \frac{1}{|I_{\mathfrak{p}}|} \sum_{\tau \in I_{\mathfrak{p}}} \operatorname{Tr}\bigl(\rho(\tau) \circ \rho(\operatorname{Frob}_{\mathfrak{p}}^k)\bigr) \\
&= \frac{1}{|I_{\mathfrak{p}}|} \sum_{\tau \in I_{\mathfrak{p}}} \operatorname{Tr}\bigl(\rho(\tau \operatorname{Frob}_{\mathfrak{p}}^k)\bigr) 
= \frac{1}{|I_{\mathfrak{p}}|} \sum_{\tau \in I_{\mathfrak{p}}} \chi_{\rho}(\tau \operatorname{Frob}_{\mathfrak{p}}^k).
\end{align*}
\end{proof}

\begin{proof}[Proof of Proposition~\ref{prop:main}]
The proof of the proposition can be obtained by carefully determining the constant terms in the calculations of the previous paper~\cite{AK}. In particular, care must be taken when handling ramified primes. Put \(I_{\mathfrak{p}} = I_{\mathfrak{p}}(L/K)\). Taking the logarithm of the statement of DRH(B), we see that DRH(A) and DRH(B) are equivalent to the validity of the following equation:
\begin{equation}\label{eq:DRH-AB}
\lim_{x \to \infty} \Biggl(m\log\log x + \sum_{N\mathfrak{p} \leq x} \sum_{\lambda_{\mathfrak{p}}} \log\bigl(1-N\mathfrak{p}^{-1/2}\lambda_{\mathfrak{p}}\bigr)^{-1} \Biggr)
= \frac{\nu(\rho)}{2}\log 2 + \log\frac{L_K^{(m)}(1/2,\rho)}{m!} - m\gamma,
\end{equation}
where \(\lambda_{\mathfrak{p}}\) runs over the eigenvalues (counted with multiplicity) of \(\rho(\operatorname{Frob}_{\mathfrak{p}})\big|_{V^{I_{\mathfrak{p}}}} : V^{I_{\mathfrak{p}}} \to V^{I_{\mathfrak{p}}}\).
By the Taylor expansion, we obtain
\[
\sum_{N\mathfrak{p} \leq x} \sum_{\lambda_{\mathfrak{p}}} \log\bigl(1-N\mathfrak{p}^{-1/2}\lambda_{\mathfrak{p}}\bigr)^{-1} = \sum_{N\mathfrak{p} \leq x} \sum_{\lambda_{\mathfrak{p}}} \sum_{k=1}^{\infty}
\frac{\lambda_{\mathfrak{p}}^k}{kN\mathfrak{p}^{k/2}}.
\]
Since Lemma~\ref{lem:Tr} yields 
\[
\sum_{\lambda_{\mathfrak{p}}} \lambda_{\mathfrak{p}}^k = \operatorname{Tr}\bigl(\rho(\operatorname{Frob}_{\mathfrak{p}}^k)\big|_{V^{I_{\mathfrak{p}}}}\bigr) =
\frac{1}{|I_{\mathfrak{p}}|} \sum_{\tau \in I_{\mathfrak{p}}} \chi_{\rho}(\tau \operatorname{Frob}_{\mathfrak{p}}^k),
\]
equation~\eqref{eq:DRH-AB} can be rewritten as follows:
\begin{equation}\label{eq:re-DRH-AB}
\sum_{N\mathfrak{p} \leq x} \sum_{\tau \in I_{\mathfrak{p}}} \sum_{k=1}^{\infty} \frac{\chi_{\rho}(\tau \operatorname{Frob}_{\mathfrak{p}}^k)}{k |I_{\mathfrak{p}}| N\mathfrak{p}^{k/2}} =
-m(\log\log x + \gamma) + \frac{\nu(\rho)}{2}\log 2 
+ \log\frac{L_K^{(m)}(1/2,\rho)}{m!} + o(1) \qquad (x\to \infty).
\end{equation}
Next, we decompose the sum over \(k\) on the left-hand side into three parts according to whether \(k=1\), \(k=2\), or \(k \geq 3\):
\begin{align*}
\sum_{N\mathfrak{p} \leq x} \sum_{\tau \in I_{\mathfrak{p}}} \sum_{k=1}^{\infty} \frac{\chi_{\rho}(\tau \operatorname{Frob}_{\mathfrak{p}}^k)}{k |I_{\mathfrak{p}}| N\mathfrak{p}^{k/2}} &=
\underbrace{
\sum_{N\mathfrak{p} \leq x} \sum_{\tau \in I_{\mathfrak{p}}} \frac{\chi_{\rho}(\tau \operatorname{Frob}_{\mathfrak{p}})}{|I_{\mathfrak{p}}| N\mathfrak{p}^{1/2}}
}_{\textcircled{1}}
+ \underbrace{
\sum_{N\mathfrak{p} \leq x} \sum_{\tau \in I_{\mathfrak{p}}} \frac{\chi_{\rho}(\tau \operatorname{Frob}_{\mathfrak{p}}^2)}{2|I_{\mathfrak{p}}| N\mathfrak{p}}
}_{\textcircled{2}}
+ \underbrace{
\sum_{N\mathfrak{p} \leq x} \sum_{\tau \in I_{\mathfrak{p}}} \sum_{k=3}^{\infty} \frac{\chi_{\rho}(\tau \operatorname{Frob}_{\mathfrak{p}}^k)}{k |I_{\mathfrak{p}}| N\mathfrak{p}^{k/2}} 
}_{\textcircled{3}}.
\end{align*}

First, the term for \(k=1\) can be rewritten as follows:
\begin{equation*}
\textcircled{1} = \sum_{\substack{N\mathfrak{p} \leq x \\ \mathfrak{p} \in S_K \setminus S_{\mathrm{ram}}}}
\frac{\chi_{\rho}(\operatorname{Frob}_{\mathfrak{p}})}{N\mathfrak{p}^{1/2}} +
c_1(\rho) + o(1) \qquad (x\to \infty),  
\end{equation*}
with
\[
c_1(\rho) = \sum_{\mathfrak{p} \in S_{\mathrm{ram}}} \sum_{\tau \in I_{\mathfrak{p}}}
\frac{\chi_{\rho}(\tau \operatorname{Frob}_{\mathfrak{p}})}{|I_{\mathfrak{p}}| N\mathfrak{p}^{1/2}}.
\]
Since there are only finitely many primes that ramify in \(L/K\), the sum \(c_1(\rho)\) is finite.
Next, we consider the sum for \(k=2\). We have
\begin{equation}\label{eq:2}
\textcircled{2} = \frac{1}{2} \sum_{\substack{N\mathfrak{p} \leq x \\ \mathfrak{p} \in S_K \setminus S_{\mathrm{ram}}}}
\frac{\chi_{\rho}(\operatorname{Frob}_{\mathfrak{p}}^2)}{N\mathfrak{p}} + \frac{1}{2} 
\sum_{ \mathfrak{p} \in S_{\mathrm{ram}}} \sum_{\tau \in I_{\mathfrak{p}}} 
\frac{\chi_{\rho}(\tau \operatorname{Frob}_{\mathfrak{p}}^2)}{|I_{\mathfrak{p}}| N\mathfrak{p}} +o(1) \quad (x\to \infty).
\end{equation}
Using the standard character identities
\begin{align*}
& \chi_{\rho}(\sigma^2) = \chi_{\operatorname{Sym}^2 \rho}(\sigma) - \chi_{\wedge^2 \rho}(\sigma) \qquad (\sigma \in G), \\
& \chi_{\operatorname{Sym}^2\rho} = \sum_{\rho' : \mathrm{irred.}} \operatorname{mult}(\rho',\operatorname{Sym}^2\rho)\chi_{\rho'}, \qquad
\chi_{\wedge^2\rho} = \sum_{\rho' : \mathrm{irred.}} \operatorname{mult}(\rho',\wedge^2\rho)\chi_{\rho'}, \\
& \nu(\rho) = \operatorname{mult}(\mathbf{1},\operatorname{Sym}^2\rho) - \operatorname{mult}(\mathbf{1},\wedge^2\rho),
\end{align*}
we obtain
\begin{align*}
\sum_{\substack{N\mathfrak{p} \leq x \\ \mathfrak{p} \in S_K \setminus S_{\mathrm{ram}}}}
\frac{\chi_{\rho}(\operatorname{Frob}_{\mathfrak{p}}^2)}{N\mathfrak{p}} 
&= \sum_{\substack{N\mathfrak{p} \leq x \\ \mathfrak{p} \in S_K \setminus S_{\mathrm{ram}}}}
\frac{1}{N\mathfrak{p}} \left(\nu(\rho) + \sum_{\substack{\rho' \neq \mathbf{1} \\ \mathrm{irred.}}} \bigl(\operatorname{mult}(\rho',\operatorname{Sym}^2\rho) - \operatorname{mult}(\rho',\wedge^2\rho)\bigr)
\chi_{\rho'}(\operatorname{Frob}_{\mathfrak{p}}) \right) \\ 
&= \nu(\rho) \left(\sum_{N\mathfrak{p} \leq x} \frac{1}{N\mathfrak{p}} - \sum_{\mathfrak{p} \in S_{\mathrm{ram}}} \frac{1}{N\mathfrak{p}} \right) \\
&\quad + \sum_{\substack{\rho' \neq \mathbf{1} \\ \mathrm{irred.}}} \left(\bigl(\operatorname{mult}(\rho',\operatorname{Sym}^2\rho) - \operatorname{mult}(\rho',\wedge^2\rho)\bigr)
\sum_{\substack{N\mathfrak{p} \leq x \\ \mathfrak{p} \in S_K \setminus S_{\mathrm{ram}}}} \frac{\chi_{\rho'}(\operatorname{Frob}_{\mathfrak{p}})}{N\mathfrak{p}} \right) +o(1) \\
&= \nu(\rho) \left(\log\log x + M_K - \sum_{\mathfrak{p} \in S_{\mathrm{ram}}} \frac{1}{N\mathfrak{p}} \right) \\
&\quad + \sum_{\substack{\rho' \neq \mathbf{1} \\ \mathrm{irred.}}} \left(\bigl(\operatorname{mult}(\rho',\operatorname{Sym}^2\rho) - \operatorname{mult}(\rho',\wedge^2\rho)\bigr)
\sum_{\substack{N\mathfrak{p} \leq x \\ \mathfrak{p} \in S_K \setminus S_{\mathrm{ram}}}} \frac{\chi_{\rho'}(\operatorname{Frob}_{\mathfrak{p}})}{N\mathfrak{p}} \right) + o(1)\\
& \hspace{12cm}  (x\to \infty),
\end{align*}
where the last equality follows from Theorem~\ref{theo:Rosen}.
Substituting the result into \eqref{eq:2}, we can rewrite the sum for \(k=2\) as follows:
\begin{equation*}
\textcircled{2} = \frac{\nu(\rho)}{2} (\log\log x + M_K) + c_2(\rho) + o(1) \qquad (x\to \infty),
\end{equation*}
with
\begin{align*}
c_2(\rho) &= -\frac{\nu(\rho)}{2} \sum_{\mathfrak{p} \in S_{\mathrm{ram}}} \frac{1}{N\mathfrak{p}} + \frac{1}{2} \sum_{\substack{\rho' \neq \mathbf{1} \\ \mathrm{irred.}}}
\left(\bigl(\operatorname{mult}(\rho',\operatorname{Sym}^2\rho) - \operatorname{mult}(\rho',\wedge^2\rho)\bigr)
\sum_{\mathfrak{p} \in S_K \setminus S_{\mathrm{ram}}}
\frac{\chi_{\rho'}(\operatorname{Frob}_{\mathfrak{p}})}{N\mathfrak{p}} \right)\\
&\quad + \frac{1}{2} \sum_{\mathfrak{p} \in S_{\mathrm{ram}}}
\sum_{\tau \in I_{\mathfrak{p}}} \frac{\chi_{\rho}(\tau \operatorname{Frob}_{\mathfrak{p}}^2)}{|I_{\mathfrak{p}}| N\mathfrak{p}}.
\end{align*}
By Lemma~\ref{lem:const}, the series defining \(c_2(\rho)\) converges.
Finally, for the sum over \(k \ge 3\), Lemma~\ref{lem:const} implies that \(\textcircled{3} = c_3(\rho) + o(1)\) as \(x \to \infty\), where the series
\[
c_3(\rho) = \sum_{\mathfrak{p} \in S_K} \sum_{\tau \in I_{\mathfrak{p}}} \sum_{k=3}^{\infty} \frac{\chi_{\rho}(\tau \operatorname{Frob}_{\mathfrak{p}}^k)}{k |I_{\mathfrak{p}}| N\mathfrak{p}^{k/2}}
\]
converges.
Setting \(c(\rho) = c_1(\rho) + c_2(\rho) + c_3(\rho)\), we obtain the claim of Proposition~\ref{prop:main} from \eqref{eq:re-DRH-AB} (together with the relation between \(M_K\) and \(c_K\)).
\end{proof}

\begin{proof}[Proof of Theorem~\ref{theo:main}]
First, we deduce (2) from (1). For any \(\sigma \in G\), the orthogonality of characters yields
\begin{align*}
[L:K] \pi_{\frac{1}{2}} (x;\sigma) &= [L:K] \sum_{\substack{\mathfrak{p} \in S_{\sigma} \\ N\mathfrak{p} \leq x}} \frac{1}{N\mathfrak{p}^{1/2}} =
|c_{\sigma}| \sum_{\substack{\mathfrak{p} \in S_K \setminus S_{\mathrm{ram}} \\ N\mathfrak{p} \leq x}} \frac{1}{N\mathfrak{p}^{1/2}} \sum_{\rho : \mathrm{irred.}}
\chi_{\rho} \left( \left(\frac{L/K}{\mathfrak{p}} \right) \right) \overline{\chi_{\rho}(\sigma)} \\
&= |c_{\sigma}| \left\{ \sum_{\substack{\mathfrak{p} \in S_K \setminus S_{\mathrm{ram}} \\ N\mathfrak{p} \leq x}} \frac{1}{N\mathfrak{p}^{1/2}} + \sum_{\substack{\rho \neq \mathbf{1} \\ \mathrm{irred.}}}
\overline{\chi_{\rho}(\sigma)} \sum_{\substack{\mathfrak{p} \in S_K \setminus S_{\mathrm{ram}} \\ N\mathfrak{p} \leq x}} \frac{\chi_{\rho} \left( \left(\frac{L/K}{\mathfrak{p}} \right) \right)}{N\mathfrak{p}^{1/2}} \right\}.
\end{align*}
Therefore, from Proposition~\ref{prop:main} and 
\[
\sum_{\substack{\mathfrak{p} \in S_K \setminus S_{\mathrm{ram}} \\ N\mathfrak{p} \leq x}} 
\frac{1}{N\mathfrak{p}^{1/2}} = \pi_{\frac{1}{2},K}(x) - \sum_{\mathfrak{p} \in S_{\mathrm{ram}}}
\frac{1}{N\mathfrak{p}^{1/2}} + o(1) \qquad (x\to \infty),
\]
we obtain (2).

Next, we deduce (1) from (2). Let \(\rho\) be a nontrivial irreducible representation of \(G\). Multiplying (2) by \(\chi_{\rho}(\sigma)\), summing over \(\sigma \in G\), and using
\[
\sum_{\sigma \in G} \chi_{\rho}(\sigma) = 0, \qquad \sum_{\sigma \in G} \chi_{\rho}(\sigma) M(\sigma) = \frac{[L:K]}{2} \, \nu(\rho), \qquad \sum_{\sigma \in G} \chi_{\rho}(\sigma) m(\sigma) = [L:K] m_{\rho},
\]
we obtain the following equation:
\begin{equation*}
\sum_{\sigma \in G} \frac{\chi_{\rho}(\sigma)}{|c_{\sigma}|} \pi_{\frac{1}{2}} (x;\sigma) = -\left( \frac{\nu(\rho)}{2} + m_{\rho} \right) \log\log x - \frac{1}{[L:K]} \sum_{\sigma \in G} 
\chi_{\rho}(\sigma) c(\sigma) + o(1) \qquad (x\to \infty).
\end{equation*}
Rewriting the left-hand side of the above equation as 
\begin{align*}
\sum_{\sigma \in G} \frac{\chi_{\rho}(\sigma)}{|c_{\sigma}|} \pi_{\frac{1}{2}} (x;\sigma) &= \frac{1}{[L:K]} 
\sum_{\sigma \in G} \chi_{\rho}(\sigma) \sum_{\substack{\mathfrak{p} \in S_K \setminus S_{\mathrm{ram}} \\ N\mathfrak{p} \leq x}} \frac{1}{N\mathfrak{p}^{1/2}} 
\sum_{\rho' : \mathrm{irred.}} \chi_{\rho'} \left( \left( \frac{L/K}{\mathfrak{p}} \right)\right) \overline{\chi_{\rho'}(\sigma)} \\
&= \sum_{\substack{\mathfrak{p} \in S_K \setminus S_{\mathrm{ram}} \\ N\mathfrak{p} \leq x}} \frac{\chi_{\rho} \left( \left( \frac{L/K}{\mathfrak{p}} \right)\right)}{N\mathfrak{p}^{1/2}} 
\end{align*}
and substituting \(c(\sigma)\) given in Theorem~\ref{theo:main}, we obtain the equation in Proposition~\ref{prop:main}, which shows that DRH(A) and DRH(B) hold for \(\rho\).
\end{proof}

\section{Numerical example}

We have made the programs and data illustrating Theorem~\ref{theo:main} for several number fields available at~\cite{AY}.
Here, we extract the example of the Gaussian number field, which corresponds to the original Chebyshev's bias.
Put \(L=\mathbb{Q}(i)\) and \(K=\mathbb{Q}\). Then
\[
G = \operatorname{Gal}(\mathbb{Q}(i)/\mathbb{Q}) = \{\sigma_1,\sigma_{-1}\},
\]
where \(\sigma_{-1}(i)=-i\). The unique nontrivial irreducible representation is one-dimensional and is identified with the primitive Dirichlet character \(\chi\) modulo \(4\).
Its Artin \(L\)-function is the Dirichlet \(L\)-function \(L(s,\chi)\), and we have \(L(1/2,\chi)=0.667691457190\ldots\).
Moreover, the values appearing in Theorem~\ref{theo:main} are given as follows:
\begin{align*}
& m_{\chi} = \operatorname{ord}_{s=1/2} L(s,\chi) = 0,\\
& \nu(\chi) = \operatorname{mult}(\mathbf{1}, \operatorname{Sym}^2 \chi) - \operatorname{mult}(\mathbf{1}, \wedge^2 \chi) = 1,\\
& M(\sigma_1) = \frac{1}{2}, \quad M(\sigma_{-1}) = -\frac{1}{2}, \quad m(\sigma_1) = m(\sigma_{-1}) = 0,\\
& c_{\mathbb{Q}} = 0.3157184520539\ldots, \\
& \alpha_{\mathbb{Q}} = \operatorname{Res}_{s=1} \zeta(s) = 1,\\
& c(\chi) = -\frac{1}{4} + \sum_{k=3}^{\infty} \sum_{p} \frac{\chi(p^k)}{kp^{k/2}} = -0.2596340772\ldots.
\end{align*}
Theorem~\ref{theo:main} asserts that DRH(A) and DRH(B) hold if and only if the following two asymptotic formulas hold:
\begin{align*}
\pi_{\frac{1}{2},\mathbb Q}(x) - 2\pi_{\frac{1}{2}}(x;\sigma_1) 
&= \sum_{p\leq x}\frac{1}{\sqrt{p}} - 2\sum_{\substack{p \leq x \\ p \equiv 1 \bmod 4}} \frac{1}{\sqrt{p}}
= \frac{1}{2}\log\log x + C(\sigma_1) + o(1), \\
\pi_{\frac{1}{2},\mathbb Q}(x) - 2\pi_{\frac{1}{2}}(x;\sigma_{-1}) 
&= \sum_{p\leq x}\frac{1}{\sqrt{p}} - 2\sum_{\substack{p \leq x \\ p \equiv 3 \bmod 4}} \frac{1}{\sqrt{p}}
= -\frac{1}{2}\log\log x + C(\sigma_{-1}) + o(1) \qquad (x\to\infty),
\end{align*}
where
\begin{align*} 
& C(\sigma_1) = \frac{1}{2}\gamma + \frac{1}{\sqrt{2}} - \frac{1}{2}\log 2 - \frac{1}{2}c_{\mathbb{Q}} - \log L(1/2,\chi) + c(\chi) = 0.63557\ldots, \\
& C(\sigma_{-1}) = -\frac{1}{2}\gamma + \frac{1}{\sqrt{2}} + \frac{1}{2}\log 2 + \frac{1}{2}c_{\mathbb{Q}} + \log L(1/2,\chi) - c(\chi) = 0.77863\ldots, \\
& C(\sigma_1) + C(\sigma_{-1}) = \sqrt{2}.
\end{align*}
This is a weighted version of Chebyshev's bias: primes congruent to \(3\) modulo \(4\) occur earlier than those congruent to \(1\) modulo \(4\).
The following table shows the values of \(\pi_{\frac{1}{2}}(x) - 2\pi_{\frac{1}{2}}(x;\sigma_1) - \frac{1}{2}\log\log x\), and 
the graph compares the numerical data of \(\pi_{\frac{1}{2}}(x) - 2\pi_{\frac{1}{2}}(x;\sigma)\) (represented by the jagged line) with the theoretical values (represented by the smooth curve) up to \(10^{15}\).

\begin{table}[H]   \centering  
\caption*{Values of $\pi_{\frac{1}{2}}(x) - 2\pi_{\frac{1}{2}}(x;\sigma_1) - \frac{1}{2}\log\log x$ .}
  \label{tab:chebyshev_bias}
\renewcommand{\arraystretch}{1.2}
\begin{tabular}{cccl}     \hline
$x$ &   $ \pi_{\frac{1}{2}}(x)-2\pi_{\frac{1}{2}} (x;\sigma_1) -\frac{1}{2} \log{\log{x}}  $\\  \hline  
$10^5$ &  $0.607110160$ \\
$10^6$  & $0.697306370$\\
$10^7$  & $0.633834715$ \\
$10^8$ & $0.618072165$ \\
$10^9$ & $0.597502570$\\
$10^{10}$ & $0.646552286$\\
$10^{11}$ & $0.637990769$\\
$10^{12}$ & $0.660303360$\\
$10^{13}$ & $0.636912376$\\
$10^{14}$ & $0.620211301$\\
$10^{15}$ & $0.634322422$ \\ \hline 
\end{tabular} \end{table}

\begin{graph}[H] \centering
 \caption*{Comparison between the numerical data and the theoretical curve.}
  \label{fig:chebyshev_bias}
\includegraphics[width=14cm]{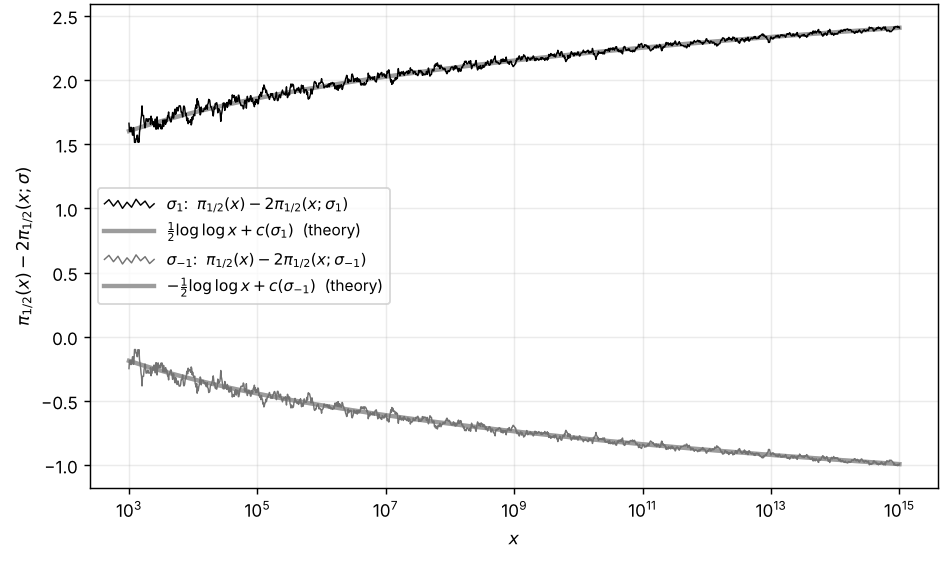}
\end{graph}

\section*{Acknowledgments}

The author is grateful to Ikuya Kaneko and Shigetoshi Yokoyama for valuable advice.

\end{document}